# $L^2$- decomposition of the second fundamental form of a hypersurface in the study of the general relativistic vacuum constraint equations

Sergey E. Stepanov[a,b,*], Irina I. Tsyganok[b]

[a] *Department of Mathematics, Russian Institute for Scientific and Technical Information of the Russian Academy of Sciences, 20, Usievicha street, 125190 Moscow, Russia*
[b] *Department of Mathematics, Finance University, 49-55, Leningradsky Prospect, 125468 Moscow, Russia*

**Abstract.** In present article, we consider a $L^2$-orthogonal decomposition of the second fundamental form of a closed spacelike hypersurface in a Lorentzian spacetime and its applications to the study of some algebraic-differential properties of the general relativistic vacuum constraint equations. For the study we will use the well-known Ahlfors Laplacian.

## 1. Introduction

In the present paper we consider a closed (i.e., compact and without boundary) Riemannian manifold $(M, g)$ of dimension $n \geq 2$ immersed in a Lorentzian manifold $\mathbb{L}_1^{n+1}$. The study of such hypersurfaces has attracted the interest of a considerable group of geometers as it is evidenced by the amount of works which was generated in the last decades. This is due not only to its mathematical interest, but also to its relevance in General Relativity. In particular, was studied restrictions on the solutions of the initial value constraints of General Relativity when the spatial hypersurface is closed.

For the study of such hypersurfaces, we will use the *Ahlfors Laplacian*. Recall that this is a second-order elliptic operator $S^*S$, defined on differential one-forms on $(M, g)$, which was considered in a number of papers [1] - [6]. The kernel of $S^*S$ is the vector space of conformal Killing one-forms. These forms are dual to conformal Killing vector fields (see [7]). In turn, any conformal Killing vector field on $(M, g)$ is a vector field whose (locally defined) flow defines conformal transformations. One of main applications of $S^*S$ is the construction of families of smooth quasi-conformal deformations of transformations $(M, g)$.

___________________________

* Corresponding author.
E-mail address: s.e.stepanov@mail.ru (S.E. Stepanov)

In the next section of the present paper, we consider the Ahlfors Laplacian in connection with the geometry of a closed Riemannian manifold. In the third section of our article, we consider a $L^2$-orthogonal decomposition of the of the second fundamental form and its traceless part of a closed spacelike hypersurface in a Lorentzian manifold, and in the fourth section we study its application to some algebraic-differential properties of the general relativistic vacuum constraint equations. For the study, we will use the Ahlfors Laplacian defined above.

We announced some of these results in our reports at the ''Seminar on current problems of geometry. Week of Lobachevsky'' (30 Nov.-2 Dec. 2023, Kazan, Russia) and "2$^{th}$ International symposium on current developments in fundamental and applied mathematics sciences" (14-17 Nov. 2023, Ankara, Türkiye).

## 2. Ahlfors Laplacian and its connection with Hodge and Sampson Laplacians

We denote by $S^p M := S^p T^* M$ the vector bundle of symmetric bilinear differential $p$-forms or, in other words, of covariant symmetric $p$-tensors $(p \geq 1)$ on $(M, g)$ and define the $L^2$ *inner product* of two covariant symmetric $p$-tensors $\varphi$ and $\varphi'$ on $(M, g)$ by the formula

$$\langle \varphi, \varphi' \rangle := \int_M g(\varphi, \varphi') dvol_g < +\infty$$

where $dvol_g$ being the volume element of $(M, g)$. Also $\delta^* : C^\infty S^1 M \to C^\infty S^2 M$ will be the first-order differential operator defined by the formula (see [8, pp. 35; 117])

$$\delta^* \theta = \frac{1}{2} L_\xi g,$$

where $L_\xi$ denotes the Lie derivatives of the vector field $\xi = \theta^\#$ (dual of the 1-form $\theta$). At the same time, we denote by the formula $\delta : C^\infty S^2 M \to C^\infty S^1 M$ the formal adjoint operator for $\delta^*$ which is called the *divergence of symmetric two-tensors* (see [8, p. 34]). In this case, we have

$$\langle \varphi, \delta^* \theta \rangle = \langle \delta \varphi, \theta \rangle \qquad (2.1)$$

for any $\varphi \in C^\infty S^2 M$ and $\theta \in C^\infty S^1 M$. In other notations we have $\delta \varphi = -div_g \varphi$.

Using (2.1), we can definite the operator $S$ from $S^1 M$ to the vector space of trace free symmetric tensors $S_0^2 M$ by the equality $S\theta := \delta^* \theta + \frac{1}{n} \delta\theta\, g$. The operator $S$ is known as the *Cauchy-Ahlfors* operator (see, for example, [3]). It's obvious that a vector field $\xi$ on $(M, g)$ is a *conformal Killing vector field* if and only if $S$ annihilates the one-form $\theta$ identified with $\xi$ by the Riemannian metric $g$, i.e., $L_\xi g = -\frac{2}{n} \delta\theta\, g$ for $\theta^\# = \xi$ (see [7]; [11, p. 46]). Therefore, the kernel of $S$ is the finite-dimensional vector space on a closed Riemannian manifold $(M, g)$.

The formal adjoint operator for $S$ is defined by the formula $S^* \omega = \delta \omega$ for an arbitrary $\omega \in C^\infty S_0^2 M$ (see [4]; [5]). Moreover, it is easy to verify the formula (see also [3]; [4] and [5]; [6, p. 337])

$$S^* S\, \theta = \frac{1}{2} \delta d\theta + \frac{n-1}{n} d\delta\theta - Ric(\xi, \cdot). \tag{2.2}$$

for the *Ahlfors operator* $S^* S: C^\infty S^1 M \to C^\infty S^1 M$ and the Ricci tensor $Ric$ of $(M, g)$. We recall that $S^* S$ is a formally self-adjoint non-negative and elliptic differential operator (see [2]; [3]; [5]; [15, pp. 24-25]). Note that $ker\, S^* S = ker\, S$ since $\langle S^* S\theta, \theta \rangle = \langle S\theta, S\theta \rangle$ for any $\theta \in C^\infty S^1 M$. Therefore, the kernel of $S^* S$ is a finite dimensional vector space of conformal Killing one-forms on a closed Riemannian manifold $(M, g)$ (see [8, p. 464]).

**Remark 1**. Note that the Ahlfors Laplacian $S^* S$ is a specific form of the second-order elliptic *Tachibana operator*. We recall that the Tachibana operator is defined on differential $p$-forms for $1 \leq p \leq n - 1$ and its kernel is the finite-dimensional vector space of conformal Killing $p$-forms. In particular, for $p = 1$ the Tachibana operator is the Ahlfors Laplacian. The geometry of the Tachibana operator is explained in our papers [9] and [10].

In (2.1), we also denote by $d$ the well known exterior derivative from the bundle $\Lambda^1 M$ of differential one-forms to the bundle $\Lambda^2 M$ of differential two-forms and by $\delta: C^\infty \Lambda^2 M \to C^\infty \Lambda^1 M$ the formal adjoint operator for $d$ which is called the *divergence of differential two-forms* (see [8, pp. 21; 54]). Then from (2.1) we obtain

$$S^*S\,\theta = \tfrac{1}{2}\triangle_H\,\theta + \tfrac{n-2}{2n}\,d\delta\theta - Ric(\xi,\cdot). \tag{2.3}$$

for the well-known *Hodge Laplacian* $\triangle_H\colon C^\infty\Lambda^1 M \to C^\infty\Lambda^1 M$, defined by the equality $\triangle_H := d\delta + \delta d$ (see [8, p. 54]). The kernel of $\triangle_H$ is the vector space of harmonic one-forms (see [8, p. 57]) and its dimension $b_1(M)$ is the first Betti number of $(M,g)$.

From (2.3) we deduce $S^*S\,\theta = \tfrac{1}{2}\triangle_H\,\theta - Ric(\xi,\cdot)$ for $n=2$. In this case, the Ahlfors operator is a Laplacian in accordance with the definition given in the well-known monograph [8, p. 52]. We also known from [12] that the Hodge Laplacian admits the following decomposition

$$\triangle_H = \triangle_S + 2Ric \tag{2.4}$$

for the *Sampson Laplacian* $\triangle_S\colon C^\infty S^1 M \to C^\infty S^1 M$ defined via the formula (see [12] and [13])

$$\triangle_S = 2\,\delta\,\delta^* - d\delta. \tag{2.5}$$

In this case, from (2.1) we deduce the following equation:

$$S^*S\,\theta = \tfrac{1}{2}\triangle_S\,\theta + \tfrac{n-2}{2n}\,d\delta\theta \tag{2.6}$$

We recall that the vector field $\xi$ is an *infinitesimal harmonic transformation* in $(M,g)$ if the local one-parameter group of infinitesimal transformations generated by the vector field $\xi$ is a local group of harmonic transformations (see [12]). Furthermore, a necessary and sufficient condition for a vector field $\xi$ on a Riemannian manifold to be infinitesimal harmonic transformation in $(M,g)$ is that $\theta \in \mathrm{Ker}\,\triangle_S$ for $\theta^\# = \xi$ (see [12]). In particular, from (2.6) we conclude that $S^*S\,\theta = \tfrac{1}{2}\triangle_S\,\theta$ for $n=2$. Therefore, any infinitesimal conformal transformation on a two-dimensional Riemannian manifold $(M,g)$ is an infinitesimal harmonic transformation (see also [12]). The converse is also true.

In addition, from (2.3) we conclude that for $n \geq 3$ and an arbitrary divergence-free one-form $\theta \in C^\infty S^1 M$ the equation $S^*S\,\theta = \tfrac{1}{2}\triangle_S\,\theta$ holds. Recall here that divergence-free vector fields on $(M,g)$ form a finite-dimensional vector space and its dimension depends only on the topology of $(M,g)$. In addition, any divergence-free infinitesimal harmonic transformation on a $n$-dimensional ($n \geq 3$) compact Riemannian manifold $(M, g)$ is an *infinitesimal isometric transformation* or, in other words, a *Killing vector field* (see also [7]; [8, p. 35]; [11, p. 43]; [12]). The converse is also true.

**Remark 2.** Considering the above, we note the following. If $\dim M \geq 3$, then from (2.1) and (2.4) we can conclude that the Ahlfors operator $S^*S: C^\infty S^1 M \to C^\infty S^1 M$ is not a Laplacian in accordance with the definition given in the well-known monograph [8, p. 52]. Namely, a second order differential operator $A$ acting on $C^\infty$-sections of a vector bundle $E$ over $(M, g)$ is a Laplace operator if its principal symbol has the form $\sigma(A)(\xi, x)\theta_x = -g(\xi, \xi)\theta_x$ for $\xi$ in $E_x$ and $\theta \in T_x^* M - \{0\}$ at an arbitrary $x \in M$. On the other hand, we deduce from (2.6) that $\sigma(S^*S)(\xi, x)\theta_x = -\frac{1}{2}g(\xi, \xi)\theta_x - \frac{n-2}{2n}\theta_x(\xi)\xi$ since from [12] we know that $\sigma(\triangle_S \theta)(\xi, x)\theta_x = -g(\xi, \xi)\theta_x$ for $\xi \in T_x M$. In addition, we have $g(-\sigma(S^*S)(\xi, x)\theta_x, \theta_x) > 0$ for $\xi \in T_x M - \{0\}$ and $\theta \in T_x^* M - \{0\}$. Therefore, $S^*S$ is an elliptic differential operator (see [8, pp. 461-463]) and $S^*S$ is not a Laplacian. However, the Ahlfors operator will still be referred to as the Ahlfors Laplacian in the present paper, due to a tradition that can be traced back to two other references [4] and [6].

## 2. $L^2$-orthogonal decomposition of the second fundamental form of a closed spacelike hypersurface in a Lorentzian spacetime

Let $\mathbb{L}_1^{n+1}$ be a Lorentzian manifold of dimension $n + 1 \geq 4$ and signature $(- + \cdots +)$ solving the Einstein field equations (see [14]; [15]; [16])

$$\overline{Ric} - \frac{1}{2} \bar{s} \bar{g} + 2\Lambda g = \kappa T$$

where we denote by $\bar{g}$ the metric tenor of $\mathbb{L}_1^{n+1}$, we also denote by $\overline{Ric}$ and $\bar{s}$ the Ricci curvature and the scalar curvature of $\bar{g}$, respectively. We employ the letter $\kappa$ for a positive constant, whose value (and physical dimensions) depends on the specific conventions one adopts. In addition, as it is rather customary in the physical literature, $T$ stands for the *stress-energy tensor* of the sources while $\Lambda$ stands for a *cosmological constant*. Any such $\mathbb{L}_1^{n+1}$ is called a *spacetime* (see, for example, [28]).

Moreover, we consider a closed spacelike hypersurface $(M, g)$ in $\mathbb{L}_1^{n+1}$ (see, for example, [8, p. 38-39]). In this case, the induced metric $g$ on $M$ from that of $\mathbb{L}_1^{n+1}$ is positive definite. The extrinsic geometry of $(M, g)$ is described by the so-called *second fundamental form* $K \in S^2 M$ (see [8, p. 38])

At the same time, it is well-known that for any $n$-dimensional ($n \geq 3$) closed Riemannian manifold $(M, g)$, the algebraic sum $\operatorname{Im} \delta^* + C^\infty M \cdot g$ is closed in $S^2 M$, and we have the decomposition

$$S^2 M = (\operatorname{Im} \delta^* + C^\infty M \cdot g) \oplus \left( \delta^{-1}(0) \cap \operatorname{trace}_g^{-1}(0) \right) \qquad (3.1)$$

where both factors are infinite dimensional and orthogonal to each other with respect to the $L^2$ inner scalar product (see [8, p. 130] and [15, pp. 24-25]). It's obvious that the second factor $\delta^{-1}(0) \cap \operatorname{trace}_g^{-1}(0)$ of (3.3) is the space of $TT$-tensors.

**Remark 3.** Here we recall that a symmetric divergence free and traceless covariant two-tensor is called $TT$-tensor (see, for instance, [17]). As a consequence of a result of Bourguignon-Ebin-Marsden (see [8, p. 132]; [18]) the space of $TT$-tensors is an infinite-dimensional vector space for any closed Riemannian manifold $(M, g)$. Such tensors are of fundamental importance in stability analysis in General Relativity (see, for instance, [19]; [15]) and in Riemannian geometry (see, for instance, [8, p. 346-347]; [18]; [20] and [21]).

Considering the above, we conclude that the second fundamental form of $K$ has the following $L^2$-orthogonal decomposition (see also formula (3.1))

$$K = \left( \frac{1}{2} L_\xi g + \lambda g \right) + \varphi^{TT} \qquad (3.2)$$

for some vector field $\xi \in C^\infty TM$, $TT$-tensor $\varphi^{TT} \in C^\infty S^2 M$ and scalar function $\lambda \in C^\infty M$. If we denote by $H := \operatorname{trace}_g K$ the *mean curvature* of $(M, g)$ (see [8, p. 38]), then applying the operator $\operatorname{trace}_g$ to both sides of (3.2), we obtain

$$H = -\delta \theta + n \lambda \qquad (3.3)$$

where $\theta^\# = \xi$. In this case, (3.2) can be rewritten in the form

$$K_0 = S\theta + \varphi^{TT} \qquad (3.4)$$

where $K_0 = K - \frac{1}{n} H\, g$ is the traceless part of the second fundamental form $K$ and $S\theta = L_\xi g + \frac{1}{n} \delta\theta\, g$ is the Cauchy-Ahlfors operator. Next, applying $S^*$ to both sides of (3.4), we obtain

$$S^*S\,\theta = \delta K_0 \qquad (3.5)$$

for the Ahlfors Laplacian $S^*S$. Therefore, $\xi = \theta^{\#}$ is a conformal Killing vector field if and only if $K_0$ is a TT-tensor. Now, we are ready to formulate our result.

**Theorem 1.** *Let $(M, g)$ be an n-dimensional closed spacelike hypersurface in a Lorentzian manifold $\mathbb{L}_1^{n+1}$ of dimension $n + 1 \geq 4$. Then there are some one-form $\theta \in C^\infty S^1 M$ and TT-tensor $\varphi^{TT} \in C^\infty S^2 M$ such that the traceless part $K_0$ of the second fundamental form $K$ of $(M, g)$ has the decomposition*

$$K_0 = S\theta + \varphi^{TT} \qquad (3.6)$$

*for the Cauchy-Ahlfors operator $S$. Furthermore, both factors in the right side of (3.6) are orthogonal to each other with respect to the $L^2$ inner product and the following equation holds*

$$S^*S\,\theta = \delta K_0 \qquad (3.7)$$

*for the Ahlfors Laplacian $S^*S$. In addition, $\xi = \theta^{\#}$ is a conformal Killing vector field if and only if $K_0$ is a TT-tensor.*

## 3. Applications to study the problem of solutions of general relativistic vacuum constraint equations

It is then a well-known fact that the triple $(M, g, K)$ solves a system of general relativistic constraint equations that takes the form (see [14]; [15]; [16, p. 44])

$$\begin{cases} s - g(K, K) + (trace_g K)^2 = 2\, T(N, N) + 2\, \Lambda; \\ div_g K - d(trace_g K) = \kappa\, T(N, \cdot). \end{cases} \qquad (4.1)$$

for a timelike unit normal vector field $N$ to $(M, g)$. This follows from the well-known Gauss–Codazzi equations for hypersurfaces, which for an arbitrary spacelike hypersurface $(M, g, K)$ isometrically immersed in a Lorentzian manifold $\mathbb{L}_1^{n+1}$.

To avoid unnecessary complications, we shall focus here on the vacuum case, i.e., we consider the field equations with no sources ($T = 0$) and take $\Lambda = 0$. In this case, the Einstein field equations can be rewritten in the form $\overline{Ric} = \frac{1}{2}\bar{s}\bar{g}$, then the field equations are also referred to as the vacuum field equations. By setting $T = 0$ and $\Lambda = 0$ in the trace-reversed field equations, the vacuum equations can be written as $Ric \equiv 0$. Then $\mathbb{L}_1^{n+1}$ is a *Ricci-flat spacetime* (see [15]). In this case, $\mathbb{L}_1^{n+1}$ is a special case of the Einstein manifold (see [8, p. 44]).

**Remark 4.** The first example of a Ricci-flat spacetime is the *Schwarzschild spacetime* that describes a static black hole. In this geometry the Ricci tensor is zero everywhere but the spacetime is most certainly not flat. The second example is a *Calabi-Yau manifold*. It is a Ricci flat manifold, but they are certainly not the same as Minkowski space $\mathbb{R}_1^{n+1}$.

The well-known problem is to construct solutions of the general relativistic vacuum constraint equations (see [15]; [16, pp. 47- 48]; [22]; [23] and etc.):

$$\begin{cases} s - g(K, K) + (trace_g K)^2 = 0; \\ div_g K - d(trace_g K) = 0. \end{cases} \quad (4.2)$$

These are the general relativistic constraint equations whatever the space-dimension $n \geq 2$. Moreover, we known the following existence theorem (see [15]).

**Theorem 2**. *Let $(M, g, K)$ be a triple satisfying the vacuum constraint equations, then there exists a spacetime $\mathbb{L}_1^{n+1}$ with metric $\bar{g}$ and an embedding $f: M \to \mathbb{L}_1^{n+1}$ such that the following assertions are true*:

1. *the spacetime $\mathbb{L}_1^{n+1}$ is Ricci-flat*;
2. *$g$ is the induced metric by $f$, namely $g = f^*(\bar{g})$*;
3. *$K$ is the second fundamental form of $f: M \to \mathbb{L}_1^{n+1}$*.

From the first equation of (3.9) we conclude that $s \geq -\frac{n-1}{n} H^2$. In particular, if $(M, g)$ is a maximal hypersurface (that is, the mean curvature of a spacelike hypersurface is zero), then its scalar curvature $s = g(K, K) \geq 0$. If in addition $s = 0$, then $(M, g)$ is a totally geodesic submanifold of $\mathbb{L}_1^{n+1}$ (that is, the second fundamental form $K$ vanishes identically).

On the other hand, the second equation of (3.9) can be rewritten in the form

$$\delta K = - dH. \qquad (4.3)$$

Combining (3.7) and (4.3), we obtain

$$S^*S\,\theta = -\frac{n-1}{n} dH.$$

Therefore, if the mean curvature $H$ of $(M,g)$ is constant, then $S^*S\,\theta = 0$ and hence $\xi = \theta^\#$ is a conformal Killing vector field. In this case, from (3.6) we obtain $K = \frac{1}{n} H\, g + \varphi^{TT}$. The converse is also true. Furthermore, if the mean curvature $H$ of $(M,g)$ is constant, then from (3.6) we obtain

$$s = g(\varphi^{TT}, \varphi^{TT}) - \frac{n-1}{n} H^2. \qquad (4.4)$$

In this case, if $(M,g)$ is a *maximal hypersurface* (that is, the mean curvature of a spacelike hypersurface is zero), then $s = g(\varphi^{TT}, \varphi^{TT}) \geq 0$. Then the following corollary holds.

**Corollary 1.** Let $\mathbb{L}_1^{n+1}$ be a vacuum spacetime and $(M,g)$ be an n-dimensional $(n \geq 3)$ closed spacelike hypersurface in $\mathbb{L}_1^{n+1}$. Moreover, let $L^2$-decomposition (3.6) holds. Then the mean curvature $H$ of $(M,g)$ is a constant if and only if $\xi = \theta^\#$ is a conformal Killing vector field and the second fundamental form $K$ of $(M,g)$ has the form $K = \frac{1}{n} H\, g + \varphi^{TT}$. In particular, if $(M,g)$ is a maximal hypersurface, then its scalar curvature $s = g(\varphi^{TT}, \varphi^{TT}) \geq 0$.

We consider Lorentzian manifolds which admit a timelike conformal vector field. Any such manifold is called a *conformally stationary spacetime* (see [29]). In turn, if $\mathbb{L}_1^{n+1}$ is a conformally stationary and Einstein spacetime with $n+1 \geq 4$, then its every closed spacelike hypersurface $(M,g)$ with constant mean curvature $H$ is *totally umbilical*, i.e., $K_0 = 0$ (see [8, p. 38]). Next, from (3.6) we obtain $S\theta + \varphi^{TT} = 0$ and hence $S^*S\theta = 0$. At the same time, the equalities $S^*S\theta = 0$ and $S\theta = 0$ are equivalent on a compact manifold $(M,g)$ and hence $\xi = \theta^\#$ is a conformal Killing vector field. In this case the second factor $\varphi^{TT}$ in (3.6) is equal to zero. In this case, from (4.4) we obtain $s = -\frac{n-1}{n} H^2 \leq 0$.

On the other hand, let $\mathbb{L}_1^{n+1}$ be a vacuum spacetime for dimension $n \geq 3$. In this case, it can be argued that if the mean curvature $H$ of $(M, g)$ is constant and the second factor $\varphi^{TT}$ in (3.6) is equal to zero, then $(M, g)$ is totally umbilical. In this case, from (4.4) we obtain $s = -\frac{n-1}{n} H^2 \leq 0$. Therefore, we can formulate the second corollary.

**Corollary 2.** *Let $\mathbb{L}_1^{n+1}$ be a vacuum spacetime for dimension $n \geq 3$ and $(M, g, K)$ be a triple satisfying the vacuum constraint equations (3.9) for $trace_g K = const$. If $\mathbb{L}_1^{n+1}$ is conformally stationary, then the second factor $\varphi^{TT}$ in the $L^2$-orthogonal decomposition (3.6) is equal to zero and the scalar curvature $s$ is a nonpositive constant. Furthermore, the vector field $\xi = \theta^{\#}$ in the $L^2$-orthogonal decomposition (3.6) is a conformal Killing vector field. Conversely, if the second factor $\varphi^{TT}$ in the $L^2$-orthogonal decomposition (3.6) is equal to zero, then $K_0 = 0$ and the scalar curvature $s$ is a nonpositive constant.*

In conclusion, we prove the following theorem.

**Theorem 3.** *Let $\mathbb{L}_1^{n+1}$ be a vacuum spacetime and $(M, g)$ be an $n$-dimensional $(n \geq 2)$ closed spacelike hypersurface in $\mathbb{L}_1^{n+1}$. Let the second fundamental form $K$ of $(M, g)$ admits $L^2$-orthogonal decomposition (3.6) for some vector field $\xi \in C^{\infty} TM$ and TT-tensor $\varphi^{TT} \in C^{\infty} S^2 M$. Then the assumption $\int_M (L_\xi H) dvol_g \geq 0$ for the mean curvature $H$ of $(M, g)$ implies that $H$ is a constant and $\xi$ is a Killing vector field.*

*Proof.* From (3.2) we deduce that

$$n \lambda = H + \delta \theta. \tag{4.5}$$

At the same time, the second equation from (4.2) can be rewritten in the form $\delta K = - dH$. Then applying the operator $\delta$ to (3.2) and using (4.5), we find

$$- dH = \delta \delta^* \theta - n \, d\lambda = \delta \delta^* \theta - (d\delta\theta + dH) =$$

$$= \Delta_S \theta - \delta \delta^* \theta - dH,$$

where $\Delta_S = 2 \delta \delta^* - d\delta$. As a result, we have $\Delta_S \theta = \delta \delta^* \theta$. In this case from (2.6) and (3.5) we deduce

$$\delta K_0 = \frac{1}{2} \delta \delta^* \theta + \frac{n-2}{2n} d \delta \theta.$$

On the other hand, we derive the following equalities:

$$\delta K_0 = \delta \left( K - \frac{1}{n} H g \right) = \delta K - \frac{1}{n} dH = -dH - \frac{1}{n} dH = -\frac{n+1}{n} dH.$$

As a result, we have

$$-\frac{n+1}{n} dH = \frac{1}{2} \delta \delta^* \theta + \frac{n-2}{2n} d \delta \theta.$$

From the last equation we obtain the integral formula

$$-\frac{n+1}{n} \int_M (L_\xi H) dvol_g = \frac{1}{2} \langle \delta^* \theta, \delta^* \theta \rangle + \frac{n-2}{2n} \langle \delta \theta, \delta \theta \rangle \geq 0.$$

Then the assumption $\int_M (L_\xi H) dvol_g \geq 0$ implies that $\delta^* \theta = 0$ and, therefore, $\delta \theta = 0$. In this case, $\xi$ is a Killing vector field and $H = constant$.

**Conclusion.** In this article, we proposed some conditions under which the mean curvature of space-like hypersurfaces will be constant. Qualitatively new in these studies is the use of the $L^2$-orthogonal decompositions of second fundamental forms of a spacelike hypersurfaces and the Ahlfors Laplacian. In turn, the importance of spacelike hypersurfaces of constant mean curvature is well known, for instance, because of the role they play in different problems in General Relativity. Namely, these hypersurfaces relevant for the study of the Einstein equation: the space of solutions of Einstein equations (see [24]) and in numerical integration schemes for Einsteins equations (see [25] and [26]). Moreover, we recall that constant mean curvature spacelike hypersurfaces have been used to prove positivity of mass (see [27]). A summary of other reasons justifying the study of these hypersurfaces in General Relativity can be found in [28].

**References**

1. Ahlfors L., Conditions for quasiconformal deformation in several variables, Contributions to Analysis. A Collection of Papers Dedicated to L. Bers, Academic Press, New York (1974) 19-25.

2. Pierzchalski A., Ricci curvature and quasiconformal deformation of a Riemannian manifold, Manuscripta Math. 66 (1989) 113-127.

3. Pierzchalski A., Gradients: the ellipticity and the elliptic boundary conditions – a jigsaw puzzle, Folia Mathematica 19 (1) (2017) 65-83.

4. Pierzchalski A., Orsted B., The Ahlfors Laplacian on a Riemannian manifold with boundary, Michigan Math. J. 43 (1996) 99-122.

5. Kozlowski W., Pierzchalski A., Natural boundary value problems for weighted form Laplacians, Ann. Scuola Norm. Sup. Pisa Cl. Sci. VII (5) (2008) 343-367.

6. Branson T., Stein-Weiss operators and ellipticity, Journal of Functional Analysis 151 (1997) 334-383.

7. Rademacher H.-B., Einstein Spaces with a Conformal Group, Results in Mathematics 56:1 (2009), 421-444.

8. Besse A.L., Einstein manifolds, Springer-Verlag, Berlin & Heidelberg, 2008.

9. Stepanov S. E., Mikesh J., The Hodge–de Rham Laplacian and Tachibana operator on a closed Riemannian manifold with curvature operator of definite sign, Izvestiya: Mathematics, London Mathematical Society 79 (2) (2015) 375–387.

10. Stepanov S.E., Mikeš J., Betti and Tachibana numbers of closed Riemannian manifolds, Differential Geometry and its Applications 31 (2013) 486-495.

11. Yano K., Integral formulas in Riemannian geometry, Marcel Dekker, New York, 1970.

12. Stepanov S.E., Shandra I.G., Geometry of infinitesimal harmonic transformations, Annals of Global Analysis and Geometry 24 (2003), 291–299.

18. Mikes J., Rovenski V., Stepanov S.E., An example of Lichnerowicz-type Laplacian, Annals of Global Analysis and Geometry 58 (2020), 19-34.

17. Gicquaud R., Ngo Q.A., A new point of view on the solutions to the Einstein constraint equations with arbitrary mean curvature and small *TT*-tensor, Class. Quant. Grav. 31:19 (2014) 195014.

19. Garattini R., Self sustained tranversable wormholes?, Class. Quant. Grav. 22: 6 (2005) 2673–2682.

18. Bourguignon, J. P. and Ebin, David G. and Marsden, Jerrold E. Sur le noyau des opérateurs pseudo-differentiels á symbole surjectif et non injectif, Comptes rendus hebdomadaires des séances de l'Académie des sciences. Séries A et B, Sciences mathématiques et Sciences physiques, 282 (1976), 867-870.

20. Boucetta M., Spectra and symmetric eigentensors of the Lichnerowicz Laplacian on $S^n$, Osaka J. Math. 46 (2009) 235–254.

14. Shoke-Bryua Y., Mathematical problems in general relativity, Uspekhi Mat. Nauk 40: 6 (1985) 3–39.

15. Carlotto A., The general relativistic constraint equations, Living Reviews in Relativity 24: 2 (2021) 1-170.

22. Chru P. T., Delay E., On mapping properties of the general relativistic constraints operator in weighted function spaces, with applications, Mémoires de la Société Mathématique de France, 2003.

23. Helmut F., Cauchy problems for the conformal vacuum field equations in General Relativity, Commun. Math. Phys. 91 (1983) 445-472.

29. Alías L.J., Romero A., Sanchez M., Spacelike hypersurfaces of constant mean curvature and Calabi-Bernstein type problems, Tôhoku Math J. 49 (1997) 337-345.

16. Avalos R., Lira J. H., The Einstein Constraint Equations, Brazil, Books in Bytes, 2021.

21. Stepanov S.E., Tsyganok I.I., Pointwise orthogonal splitting of the space of *TT*-tensors, Differential geometry of manifolds of figures, 54:2 (2023), 45-53.

25. Eardley D., Smarr L., Time functions in numerical relativity: Marginally bound dust collapse. Phys. Rev. D, 19 (1979), 2239–2259.

24. Fischer A., Marsden J., Moncrief V., The structure of the space of solutions of Einstein's equations. I. One Killing field. Ann. Inst. H. Poincare, 33 (1980), 147–194.

26. Piran T., Problems and solutions in numerical relativity. Ann. NY Acad., 375 (1981), 1–14.

27. Schoen R., Yau S. T., Proof of the positive mass theorem. I. Commun. Math. Phys., 65 (1979), 45–76.

28. Marsden, J. E. and Tipler, F. J.: Maximal hypersurfaces and foliations of constant mean curvature in General Relativity, Phys. Rep., 66 (1980), 109–139.